\documentclass[12pt]{amsart}
\usepackage{amsmath,amsthm,latexsym,amscd,amsbsy,amssymb,url}
\usepackage[all]{xypic}
 \relax

\chardef\bslash=`\\ 

\makeatletter
\def\verbatim{\interlinepenalty\@M \@verbatim
  \leftskip\@totalleftmargin\advance\leftskip2pc
  \frenchspacing\@vobeyspaces \@xverbatim}
\makeatother
\newtheorem{thm}{Theorem}[section]
\newtheorem{cor}[thm]{Corollary}
\newtheorem{lem}[thm]{Lemma}

\newtheorem{que}[thm]{Question}

\begin{document}


\title
{On Homeomorphisms of Function Spaces Over Products with Compacta}
\author{Raushan  Buzyakova}
\email{Raushan\_Buzyakova@yahoo.com}



\begin{abstract}
We show that if $C_p(X\times Z)$ is homeomorphic to $C_p(Y\times Z)$, where $Z$ is compact,  and $X$ and $Y$ are of countable netweight, then 
$C_p(X\times M)$ is homeomorphic to $C_p(Y\times M)$ for some metric compactum $M$.
\end{abstract}

\maketitle
\markboth{R. Buzyakova}{On Homeomorphism of Function Spaces Over Products with Compacta}
{ }

\section{Introduction}\label{S:introduction}
\par\bigskip\noindent
One of the central problems in $C_p$-theory is identifying properties of $X$ and $Y$ that guarantee the existence of a homeomorphism between $C_p(X)$ and $C_p(Y)$. Let us treat this pair of function spaces as $C_p(X\times \{0\})$ and $C_p(Y\times \{0\})$. With this alternative point of view, it is justifiable  to search for conditions that lead to a homeomorphism between $C_p(X\times M)$ and $C_p(Y\times M)$ for some relatively small and not too loose space $M$, such as a metric compactum, for example. 

\par\bigskip\noindent
It was proved in \cite{Buz} that if $C_p(X\times (\omega_1 +1))$ is homeomorphic to $C_p(Y\times (\omega_1+1))$ for metric compacta $X$ and $Y$, then $C_p(X\times \lambda)$ is homeomorphic to $C_p(Y\times \lambda)$ for an isolated $\lambda<\omega_1$. In this paper, we will generalize the mentioned result. Namely,  we prove that  if $C_p(X\times Z)$ is homeomorphic to $C_p(Y\times Z)$ , where $Z$ is compact, and $X$ and $Y$ are of countable netweight, then 
$C_p(X\times M)$ is homeomorphic to $C_p(Y\times M)$ for some metric compactum $M$. The space $M$ in the statement is a continuous image of $Z$. Since any continuous metric image of $(\omega_1+1)$ is  homeomorphic to an isolated countable ordinal, the earlier result follows from the result to be proved in this paper. Note that if we do not place any restriction on Z in the hypothesis of our main result, then the conclusion is no longer true. Indeed, let $X$ be the space of reals, $Y$ be the space of irrationals, and $Z= X^\omega\times Y^\omega$. It is clear that $C_p(X\times Z)$ is homeomorphic to $C_p(Y\times Z)$. However, due to Okunev's theorem (see, for example, \cite{Oku} or  \cite[Corollary III.2.12]{Arh}),  $C_p(X\times M)$ is not homeomorphic to $C_p(Y\times M)$ for any compactum $M$ since $X$ is $\sigma$-compact while $Y$ is not.  In notation and terminology of general topological nature, we follow \cite{Eng}. For basic facts and terminology related to function spaces, we refer to \cite{Arh}.  All spaces under consideration are Tychonoff. The identity map on $X$ will be denoted by $id_X$. Recall that  a space $X$ has a countable netweight if it has a countable collection $\mathcal S$ of closed sets such that for any $x\in X$ and any open $U$ containing $x$, there exists $S\in \mathcal S$ such that $x\in S\subset U$.

\section{Main Result}\label{S:mainresult}

\par\bigskip\noindent
For clarity of exposition, let us agree on  notations for some classical structures tailored to our arguments and  to be used for the duration of this section.

\par\bigskip\noindent
{\bf Definition of $f_x$:} {\it Given $x\in X$ and $f\in C_p(X\times Z)$, the map $f_{x,Z}:Z\to \mathbb R$ is defined by letting $f_{x,Z}(z) = f(x,z)$. When  it is clear what $Z$ is  under consideration, we will  write $f_x$.
}

\par\bigskip\noindent
{\bf Definition of $f_{F,D}$:} {\it
Given $D\subset X$ and $F\subset C_p(X\times Z)$, put  $f_{D,F}=\Delta \{f_d: f\in F, d\in D\}: Z\to \mathbb R^{|D||F|}$.
}

\par\bigskip\noindent
{\bf Definition of $\Phi_{h,X}$:} {\it
Given $X$ and a continuous surjection $h: Z\to M$, define  $\Phi_{h,X}: C_p(X\times M) \to C_p(X\times Z)$ by letting $\Phi_{h,X}(f)=f\circ (id_X \times h)$. If it is clear what $X$ is under consideration, we will write $\Phi_h$.
}

\par\bigskip\noindent
Note that in the terminology of \cite[chapter 0 section 4]{Arh}, $\Phi_h$ is the restriction of $(id_X \times h)^\#$ to $C_p(X\times M)$. 

\par\bigskip\noindent
We will use the following classical facts customized to our scenarios.
\par\bigskip\noindent
{\bf Facts.}
{\it
\begin{enumerate}
	
	\item If $cf(\tau >\omega)$ and $f$ is a continuous surjection of $\tau$ onto a separable metric space $M$, then $id_X\times f: X\times \tau\to X\times M$ is a quotient map for any space $X$.

\par\smallskip\noindent
{\rm
To see why this is true, observe that there exists a separable compact subspace $C\subset \tau$ such that $(id_X\times f) (X\times C)=(id_X\times f) (X\times \tau)$. Since $id_X\times (f|_C)$ is quotient, so is $id_X\times f$ (\cite[Corollary 2.4.6]{Eng}). 
}
	\item  If $f:Z\to M$ is a perfect map, then $id_X\times f: X\times Z\to X\times M$ is a quotient map for any space $X$.

\par\smallskip\noindent
{\rm
Note that any identity map is perfect too. Hence, $id_X\times f$ is perfect. Since any perfect map is quotient, the statement follows.
}
	\item {\rm (\cite[Corollary 0.4.8, parts (1) and (2)]{Arh})} If $id_X\times h: X\times Z \to X\times M$ is a quotient map, then $\Phi_h$ is a closed embedding.
\end{enumerate}
}

\par\bigskip\noindent
The following technical folklore statement is given with a detailed proof due to its significance in the argument of the main result.
\begin{lem}\label{lem:factor}
Let $X$ be separable, $Z$ a space,  and  $f:X\times Z\to \mathbb R^\omega$  a perfect map. Then there exist continuous $F_1: Z\to \mathbb R^\omega$ and continuous $F_2: X\times F_1(Z)\to \mathbb R^\omega$ such that $f=F_2\circ (id_X\times F_1)$.
\end{lem}
\begin{proof}
Fix countable dense $A\subset X$ and let $F_1=\Delta\{f_{a,Z}: a\in A\}:  Z\to \mathbb R^\omega$. 

\par\medskip\noindent
Define $F_2$ by letting $F_2(x,p) = f(x,z)$, where $z\in F_1^{-1}(p)$. To show that $F_2(x,p)$ does not depend on our choice of $z$, pick distinct $z_1$ and $z_2$ in $F_1^{-1}(p)$. Since $F_1(z_1)=F_1(z_2)$, we conclude that $f(a,z_1)=f(a,z_2)$ for all $a\in A$. Since $A$ is dense in $X$, $f(x,z_1)=f(x,z_2)$.

\par\medskip\noindent
Let us show that $f=F_2\circ (id_X \times F_1)$. Fix $(x,z)\in X\times Z$. We have $F_2( (id_X \times F_1)(x,z))=F_2( x,F_1(z))$. Since $z\in F_1^{-1}(z)$, by the definition of $F_2$, we have $F_2( x, F_1(z))=f(x,z)$. 

\par\medskip\noindent
The function $F_1$ is continuous as a diagonal product of continuous functions. Moreover, $F_1$ is perfect since $f_{a,Z}$ is a restriction of the perfect map $f$ to the closed subspace $\{a\}\times Z$. 

\par\medskip\noindent
To show that $F_2$ is continuous, fix an open $V\subset \mathbb R^\omega$. Put $W=F^{-1}_2(V)$ and $U=f^{-1}(V)$. We need to show that $W$ is open, Since $f$ is continuous, $U$ is open. Since $id_X\times F_1$ is quotient (due to being perfect), it suffices to show that $U=(id_X\times F_1)^{-1}(W)$.  To show that $U\subset (id_X\times F_1)^{-1}(W)$, fix $(x,z)\in U$. Then $f(x,z)\in V$. By the definition of $F_2$, we conclude that $F_2((id_X\times F_1)(x,z)) = f(x,z)\in  V$. Hence, $(x,z)\in (id_X\times F_1)^{-1}(W)$. Next,  to show that $U\supset (id_X\times F_1)^{-1}(W)$, fix $(x,p)\in W$. Then, by the definition of $F_2$, we have  $f(x,z)=F_2(x,p)\in V$ for every $z\in F_1^{-1}(p)$. Therefore, $ (id_X\times F_1)^{-1}(x,p)\subset U$.
\end{proof}

\par\bigskip\noindent
Note that in the proof of Lemma \ref{lem:factor}, we constructed $F_1$ by selecting a random countable subset of $X$. Therefore, the statement of the  lemma can be written in a more constructive way, which suits us better for future references, as follows.

\par\bigskip\noindent
\begin{lem}\label{lem:factor1}
Let $X$ be separable, A a countable dense subspace of $X$, $Z$ a space,  $f:X\times Z\to \mathbb R^\omega$  a perfect map, and $G_1=\Delta\{f_{a,Z}: a\in A\}:  Z\to \mathbb R^\omega$. Then there exists  continuous $G_2: X\times G_1(Z)\to \mathbb R^\omega$ such that $f=G_2\circ (id_X\times G_1)$. 
\end{lem}

\par\bigskip\noindent
We are now ready for our main result. During the proof we will use a classical concept concerning diagonal maps. Let $I$ be an index set and  let $f_i$ be a map from $X$ to $\mathbb R$ for each $i\in I$. Given $J\subset I$, the map $p$  defined by letting 
$p(\langle f_i(x): i\in I\rangle )= \langle f_i(x): i\in J\rangle$ will be called the {\it natural projection} of  $(\Delta_{i\in I}f_i)(X)$ onto $(\Delta_{i\in J}f_i)(X)$. Participating maps and indices will be clear from the context.
\par\bigskip\noindent
\begin{thm}\label{thm:main}
Let $X$ and $Y$ have countable netweight and let $Z$ admit a perfect map onto a separable metric space. If $C_p(X\times Z)$ is homeomorphic to $C_p(Y\times Z)$, then $C_p(X\times M)$ is homeomorphic to $C_p(Y\times M)$ for some separable metric space $M$ that is a perfect image of $Z$.
\end{thm}
\begin{proof}
Fix a homeomorphism $\phi: C_p(X\times Z)\to C_p(Y\times Z)$ and countable dense subsets $A$ and $B$ of $X$ and $Y$, respectively. We will construct $M$ recursively.

\par\bigskip\noindent
{\it Step $0$.} Let $h_0: Z\to \mathbb R^\omega$ be any perfect map. Put $F_0 = \Phi_{h_0}(C_p(X\times h_0(Z)))\subset C_p(X\times Z)$.

\par\bigskip\noindent
{\it Assumption.}  Assume that for each $i<n$,   we have defined $F_i, h_i$ with the following properties:

\begin{description}
	\item[\rm Property 1] $F_i$ has countable netweight.	
	\item[\rm Property 2] For each $i<n-1$, $\phi (F_i) \subset F_{i+1}$ if $i$ is even and  $\phi^{-1} (F_i) \subset F_{i+1}$ if $i$ is odd.  	
\end{description}
\par\medskip\noindent
Before we start our inductive case, let us show that Property 2 implies the following property.
\begin{description}
	\item[\rm Property 3] For each $i < n-2$, $F_i\subset F_{i+2}$.
\end{description}
To see why Property 3 can be derived from Property 2, let $i$ be even. By Property 2 for even indices, $\phi (F_i)\subset F_{i+1}$. Therefore, $F_i\subset \phi^{-1}(F_{i+1})$. Since $i+1$ is odd, by Property 2 for odd indices, $\phi^{-1}(F_{i+1})\subset F_{i+2}$, which proves Property 3.

\par\bigskip\noindent
{\it Step $n$.} We have two cases.
\begin{description}
	\item[\underline {\rm Case "$n$ is odd"}] Since $F_{n-1}$ has countable netweight, we can fix a countable dense $D\subset \phi (F_{n-1})$. 
Put $h_n=\Delta\{f_b: f\in D, b\in B\}: Z\to \mathbb R^\omega$. Let $h = h_0\Delta .... \Delta h_n$ and $F_n = \Phi_h (C_p(Y\times h(Z)))\subset C_p(Y\times Z)$. 

\par\medskip\noindent
Since $Y\times h(Z)$ has countable netweight, so does $F_n$. Hence, Property 1 holds. 
To verify Property 2, we need to show that $\phi (F_{n-1}) \subset F_n$. Since $h$ is perfect, by Facts 2 and 3, $F_n$ is closed. Since $D$ is dense in $\phi (F_{n-1})$, it suffices to show that $D\subset F_n$. Fix $f\in D$. Put $G_1=\Delta \{f_{b, Z}: b\in B\}: Z\to \mathbb R^\omega$. By Lemma \ref{lem:factor1}, we can fix $G_2:Y\times G_1(Z)\to \mathbb R^\omega$ such that $f=G_2\circ (id_Y\times G_1)$. Note that the maps of $G_1$ participate in the definition of $h_n$, and therefore, in the definition of $h$. We can, therefore,
project $h(Z)$ onto $G_1(Z)$ by virtue of the natural projection $p$. We now have the following:

$$
f=G_2\circ (id_Y\times G_1) = G_2\circ [(id_Y\times p)\circ (id_Y\times h)]
$$
$$=$$
$$
 [G_2\circ (id_Y\times p)]\circ (id_Y\times h).
$$
The map $g = G_2\circ (id_Y\times p)$ is in $C_p(Y\times h(Z))$. Therefore, $f=\Phi_h(g)\in F_n$.

	\item[\underline{\rm Case "$n$ is even"}] Since $F_{n-1}$ has a countable netweight, we can fix a countable dense $D\subset \phi^{-1} (F_{n-1})$. 
Put $h_n=\Delta\{f_a: f\in D, a\in A\}: Z\to \mathbb R^\omega$. Let $h = h_0\Delta .... \Delta h_n$ and $F_n = \Phi_h (C_p(X\times h(Z)))$.
\end{description}

\par\bigskip\noindent
Put $\mu = \Delta \{h_n:n\in \omega\}$,  $M=\mu (Z)$, $F_X=\Phi_{\mu , X} (C_p(X\times M))$, and $F_Y=\Phi_{\mu , Y} (C_p(Y\times M))$. Since $F_X$ and $F_Y$ are homeomorphic to $C_p(X\times M)$ and $C_p(Y\times M)$, respectively, it remains to show that $F_X$ is homeomorphic to $F_Y$.

\par\bigskip\noindent
{\it Claim 1.} $F_X$ and $F_Y$ are closed in their superspaces.
\par\smallskip\noindent
To prove the claim observe that $\mu $ is a perfect map since $h_0$ is perfect. Therefore, $id_X\times \mu$ is a quotient map. By \cite[part (2) of  Corollary 0.4.8]{Arh}, $\Phi_{\mu , X} (C_p(X\times M))$ is closed in $C_p(X\times Z)$.

\par\bigskip\noindent
{\it Claim 2.} $\bigcup_{n\in \omega} F_{2n}$ is a dense subset of  $F_X$ and $\bigcup_{n\in \omega} F_{2n+1}$ is a dense subset of  $F_Y$.
\par\smallskip\noindent
To prove the claim, let us first show that $\bigcup_{n\in \omega} F_{2n}$ is a subset of $F_X$. Fix $f\in F_{2n}$ and put $h=h_0\Delta ... \Delta h_{2n}$. Then $f= g\circ (id_X\times h)$ for some $g\in C_p(X\times h(Z))$. Let $p$ be the natural projection of $\mu (Z)$ onto  $h(Z)$. Then, $g\circ (id_X\times p)\in C_p(X\times M)$ and $\Phi_{\mu , X} (g\circ (id_X\times p))=f$. Hence, $f\in F_X$.

To show that $\bigcup_{n\in \omega} F_{2n}$ is  dense in $F_X$, fix $f\in F_X$. Let $U=U(p_1,...,p_k, B_1,...,B_k)$  be a basic open neighborhood of $f$. There exists $2m$ such that  for any $i<j\leq k$,  the map $g_1=id_X\times (h_0\Delta ...\Delta h_{2m})$ distinguishes  $p_i$ and $p_j$  if and only if  $id_X\times \mu$ does . Clearly,  there exists a continuous  $g_2: g_1(X\times Z)\to \mathbb R$ such that $g_2(g_1(p_1))\in B_1,..., g_2(g_1(p_k))\in B_k$. Hence, $g_2\circ g_1$ is in $F_X\cap U$.
\par\bigskip\noindent
{\it Claim 3.} $\phi (F_X) = F_Y$.
\par\smallskip\noindent
By Properties 2 and 3 of the induction assumptions we have  $\phi (\bigcup_n F_{2n}) = \bigcup_n F_{2n+1}$. Next, by Claims 1 and 2 and the fact that $\phi$ is a homeomorphism, $\phi (F_X)=F_Y$. The claim is proved.

\par\bigskip\noindent
By Claim 3, $F_X$ and $F_Y$ are homeomorphic. Recall that   $F_X=\Phi_{\mu , X} (C_p(X\times M))$, $F_Y=\Phi_{\mu , Y} (C_p(Y\times M))$, and  the maps $\Phi_{\mu , X}$ and $\Phi_{\mu , X}$ are embeddings. Hence, $C_p(X\times M)$ is homeomorphic to $C_p(Y\times M) $.
\end{proof}

\par\bigskip\noindent
The following is an  immediate corollary to Theorem \ref{thm:main}. 
\par\bigskip\noindent
\begin{thm}\label{thm:compact}
Let $X$ and $Y$ have countable netweight and let $Z$ be compact. If $C_p(X\times Z)$ is homeomorphic to $C_p(Y\times Z)$, then $C_p(X\times M)$ is homeomorphic to $C_p(Y\times M)$ for some  metric compact space $M$.
\end{thm}

\par\bigskip\noindent
In the proof of Theorem \ref{thm:main}, we need the existence of a perfect map of $Z$ onto a separable metric space only to show that $id_X\times \mu$ is quotient, where $\mu$ is a specially constructed continuous surjection of  $Z$ onto a subspace of $\mathbb R^\omega$. The same conclusion also holds if we replace $Z$ with an ordinal $\tau$ of uncountable cofinality (see Fact 2). Therefore, our proof generalizes an earlier result of the author \cite{Buz} as follows:

\par\bigskip\noindent
\begin{cor}\label{cor:ccordinal}
Let $X$ and $Y$ have countable netweight. If $C_p(X\times \tau)$ is homeomorphic to $C_p(Y\times \tau )$, where $\tau$ is a countably compact ordinal, then 
$C_p(X\times \lambda)$ is homeomorphic to $C_p(Y\times \lambda)$ for some countable isolated ordinal $\lambda$.
\end{cor}

\par\bigskip\noindent
Note that in Corollary \ref{cor:ccordinal}, $\lambda$ plays the same  role as $M$ in  Theorem \ref{thm:main}, which is a continuous metric image of $Z$. This level of specificity  in Corollary \ref{cor:ccordinal} is achieved due to the fact that any continuous metric image of a countably compact ordinal is homeomorphic to a countable isolated ordinal.

\par\bigskip\noindent
Our main result and corollaries prompts the following open-ended questions.
\par\bigskip\noindent
\begin{que}
Let $X$ and $Y$ have countable netweights. Suppose that $C_p(X\times Z)$ is homeomorphic to $C_p(Y\times Z)$ for some $Z$ (not necessarily compact). What additional conditions on $X,Y$ guarantee that 
$C_p(X\times M)$ homeomorphic to $C_p(Y\times M)$ for some $M\subset \mathbb R^\omega$?
\end{que}

\par\bigskip\noindent
\begin{que}
Let $X$ and $Y$ have countable netweights. What additional conditions on $X,Y$ guarantee that  
$C_p(X\times M)$ is homeomorphic to $C_p(Y\times M)$ for some $M\subset \mathbb R^\omega$? 
\end{que}
\par\bigskip\noindent
In \cite{Buz}, it was observed that it is a corollary to  Okunev's result in \cite{Oku2}  that the existence of a homeomorphism between $C_p(X\times \omega_1)$ and $C_p(Y\times \omega_1)$ guarantees the existence of a homeomorphism between $C_p(X\times (\omega_1+1))$ and $C_p(Y\times (\omega_1+1))$. It is not clear, however, if the converse holds too. In general, the following question might be of interest.

\par\bigskip\noindent
\begin{que}
Let $Z$ be a pseudocompact space and $C_p(X\times \beta Z)$ is homeomorphic to $C_p(Y\times \beta Z)$. Is then $C_p(X\times Z)$ homeomorphic to $C_p(Y\times Z)$ (perhaps, under some additional conditions)?
\end{que}

\par

\end{document}